\documentclass{amsart}
\usepackage{graphicx}
\usepackage{float}
\usepackage[arc,all]{xy}
\usepackage{amsthm}
\newtheorem{thm}{Theorem}[section]
\newtheorem{cor}[thm]{Corollary}

\theoremstyle{definition}

\theoremstyle{remark}

\numberwithin{equation}{section}

\begin{document}

\title[the structure  of  nilpotent Lie algebras]{the structure  of nilpotent  Lie algebras of class two with derived subalgebra of dimension two }
\author[F. Johari]{Farangis Johari}
\address{Departamento de Matem\'atica, Instituto de Ci\^encias Exatas, Universidade  Federal de Minas Gerais, Av. Ant\^onio Carlos 6627, Belo Horizonte, MG, Brazil}
\email{farangisjohari@umfg.br, farangisjohari85@gmail.com} 
\author[A. Shamsaki]{Afsaneh Shamsaki}
\address{School of Mathematics and Computer Science\\
Damghan University, Damghan, Iran}
\email{Shamsaki.Afsaneh@yahoo.com}
\author[P. Niroomand]{Peyman Niroomand}
\address{School of Mathematics and Computer Science\\
Damghan University, Damghan, Iran}
\email{niroomand@du.ac.ir, p$\_$niroomand@yahoo.com}

%\thanks{0}%
%\subjclass{0}%


%\date{0}%
%\dedicatory{0}%
%\commby{0}%
% ----------------------------------------------------------------
\begin{abstract}
In  this paper, the structure  of  all finite-dimensional nilpotent Lie algebras of class two  with  derived subalgebra of dimension two  over an arbitrary field $ \mathbb{F} $  is determined.  Furthermore, we give  the structure of  the Schur multiplier of such Lie algebras.
\end{abstract}
\maketitle
% ----------------------------------------------------------------
\section{ introduction }
The classification of nilpotent Lie algebras  is a classical problem. The answer to this problem is not easy over an arbitrary field without putting any condition on a nilpotent Lie algebra $ L. $  Sometimes it is possible  to characterize the structure of $ L $ when conditions are put on $ L. $ In \cite{N1}, the structure of a nilpotent Lie algebra $ L $ was given when $ \dim L^{2}=1. $ When $ \dim L^{2}=2 $ and $ L $ is nilpotent of class $ 3, $ the structure of $ L $ was given in \cite{N3}, and when $ \dim L^{2}=3 $ and $ L $ is nilpotent of  class $ 4, $ the structure of $ L $ was obtained in \cite{N4}. 
Also, the structure of  capable finite-dimensional nilpotent  Lie algebras of class two with  derived subalgebra of dimension two   was characterized  in \cite{N2}.\\
The purpose of this paper is  to complete the classification of  nilpotent Lie algebras of class two with  derived subalgebra of dimension two  over an arbitrary field, using the method  from Hardy and Stitzinger in \cite{S,H}. Moreover, we give the structure of the Schur multiplier of these Lie algebras. \\
 It is known that from \cite[Proposition 2.4 ]{N2}, every finite-dimensional nilpotent Lie algebra $ L $ of  class two with $ \dim L^{2}=2 $ is descendants of smaller Lie algebras $ H $ and $ A $ where $ L=H\oplus A $ for a   nilpotent stem Lie algebra $H$ of  class two with derived subalgebra of dimension two  and an abelian Lie algebra $ A. $ Hence the classification of a nilpotent Lie algebra $L$
 depends only on the structure of the  nilpotent stem Lie algebra $H$.\\ 
A Lie algebra $ H $ is called generalized Heisenberg  of rank $ n $ if $ H^{2}=Z(H) $ and $ \dim H^{2}=n $ and 
 a  Lie algebra $H$ is  stem when $ Z(H) $ is a subset of $ H^{2}. $ Clearly, nilpotent stem Lie algebras of class two with  derived subalgebra of dimension two are equal to  generalized Heisenberg Lie algebras of rank $ 2. $ \\
In \cite[Theorems 2.7 and 3.1]{N3}, the structure of   $ n $-dimensional generalized Heisenberg Lie algebras of rank $ 2 $ over an arbitrary field for all $ n $ such that $ n\leq 7 $ was obtained as follow:
\[L_{5,8}  = \langle x_{1},\dots,x_{5}\mid [x_{1},x_{2}]=x_{4},[x_{1},x_{3}]=x_{5}\rangle , \]
\[L_{6,22}(\varepsilon)=\langle x_{1},\dots,x_{6}\mid [x_{1},x_{2}]=x_{5},[x_{1},x_{3}]=x_{6}, [x_{2},x_{4}]=\varepsilon x_{6},[x_{3},x_{4}]=x_{5}\rangle , ~~\varepsilon \in \mathbb{F},\]
\[L^{(2)}_{6,7}(\eta)=\langle x_{1},\dots,x_{6}\mid [x_{1},x_{2}]=x_{5},[x_{1},x_{3}]=x_{6}, [x_{2}, x_{4}]=\eta x_{6}, [x_{3}, x_{4}]=x_{5}+x_{6} \rangle\]
where $\eta \in {\lbrace 0, \omega \rbrace}$,
\[L_{1}=\langle x_{1}, \dots, x_{7}| [x_{1},x_{2}]=x_{6}, [x_{1},x_{4}]=x_{7}=[x_{3},x_{5}] \rangle ,\]
\[L_{2}=\langle x_{1}, \dots, x_{7}| [x_{1},x_{2}]=x_{6}=[x_{3},x_{4}], [x_{1},x_{5}]=x_{7}=[x_{2},x_{3}] \rangle .\]
\section{Main results}
 In this section, we are going to obtain the structure of all $ n $-dimensional   nilpotent  Lie algebras of class two with   derived subalgebra of dimension two for all $ n $ such that  $ n\geq 8. $  Throughout this section, $ H(m) $ denotes the Heisenberg Lie algebra of dimension $ 2m+1, $ that is, a Lie algebra $ L $ with $ L^{2}=Z(L) $ and $ \dim L^{2}=1 $ and $ A(n) $ is an $ n $-dimensional  abelian Lie algebra.\\
 First, in the following theorem, we obtain the structure  of all $ n $-dimensional generalized Heisenberg Lie algebras  of rank $ 2 $ over an arbitrary field $ \mathbb{F} $  for all $ n $ such that $ n\geq 8. $ 
\begin{thm}\label{T1}
Let $L$ be an $ n $-dimensional  generalized Heisenberg Lie algebra  of rank $ 2 $  over an arbitrary field $ \mathbb{F} $  for all  $ n $ such that $ n\geq 8. $  Then $L $ is isomorphic to one of the following Lie algebras.
 \[H_{1}=\langle x_{i}, y_{i}, z, z_{1}\mid [x_{i}, y_{i}]=z, [x_{1}, x_{2}]=z_{1}, ~~1\leq i\leq m\rangle ~~\text{for all}~~ m ~~\text{such that }~~m\geq 3.\]
 \[H_{2}=\langle x_{i}, y_{i}, z, z_{1}, q\mid [x_{i}, y_{i}]=z, [q, x_{1}]=z_{1}, ~~1\leq i\leq m\rangle ~~\text{for all}~~ m ~~\text{such that }~~ m\geq 3. \]
\[H_{3}=\langle x_{i}, y_{i}, z, z_{1}, q_{j}\mid [x_{i}, y_{i}]=z, [x_{1}, x_{2}]=[q_{j}, x_{j+2}]=z_{1}, ~~1\leq i\leq m, 1\leq j \leq k\rangle \]\[~  \text{for all}~~ k,m ~~\text{such that}~~2\leq k\leq m-2 ~\text{and}~ m\geq 4.\]
\[ H_{4}=\langle x_{i}, y_{i}, q_{j}, z, z_{1}\mid [x_{i}, y_{i}]=z, [q_{j},x_{j}]=z_{1}, 1\leq i\leq m, 1\leq j\leq k \rangle \]\[ \text{for all}~~ k,m ~~\text{such that}~~2\leq k \leq m~~ and~~ m\geq 2.\]
\[  H_{5}=\langle x_{i}, y_{i}, q_{j}, q_{j}^{\prime}, z, z_{1}\mid [x_{i}, y_{i}]=z, [x_{1}, x_{2}]=[q_{j},q^{\prime}_{j}]=z_{1}, 1\leq i\leq m, 1\leq j\leq k_{1} \rangle \]\[  \text{for all}~~ k_{1},m ~~\text{such that}~~m\geq 2~\text{and}~k_{1}\geq 2.\]
\[  H_{6}=\langle x_{i}, y_{i}, q_{j}, q_{j}^{\prime}, p_{t},z, z_{1} \mid [x_{i}, y_{i}]=z, [x_{1}, x_{2}]=[q_{j},q^{\prime}_{j}]=z_{1}, [p_{s},x_{s+2}]=z, \]\[ 1\leq i\leq m,1\leq j\leq k_{1}, 1\leq t\leq r \rangle~\]\[\text{for all}~r,m,k_{1}\text{such that}~~1\leq  r \leq m-2, m\geq 4 ~\text{and}~k_{1}\geq 2. \]
\[ H_{7}=\langle x_{i}, y_{i}, q_{j}, q_{j}^{\prime}, p_{t},z, z_{1}\mid [x_{i}, y_{i}]=z,[q_{j},q^{\prime}_{j}]=[p_{t},x_{t}]=z_{1}, \]\[1\leq i\leq m, 1\leq j\leq k_{1}, 1\leq t\leq r\rangle\]\[~~\text{for all}~ r,m,k_{1}~~\text{such that}~~1\leq r \leq m,~~ m\geq 2~~\text{and}~~k_{1}\geq 2 .\]
\[  H_{8}=\langle x_{i}, y_{i}, q_{j}, q_{j}^{\prime}, z, z_{1}\mid [x_{i}, y_{i}]=z, [q_{j},q^{\prime}_{j}]=z_{1}, 1\leq i\leq m, 1\leq j\leq k_{1}\rangle\]\[\cong H(m)\oplus H(k_{1})~~\text{for all}~k_{1},m~~\text{such that}~~ k_{1}\geq 2 ~~\text{and} ~~m\geq 1.\]
\begin{proof}
Let  $N$ be a one-dimensional central ideal of $ L $ contained in $L^2.$  Since $ L $ is of nilpotency class two, $  L^{2}=Z(L)\cong A(2) $ and $\dim L/N\geq 7.$ Since $ \dim (L/N)^{2}=1, $  \cite[Lemma 3.3]{N1} implies  $ L/N \cong H(m)\oplus A(k)  $ for all $ m,k $ such that $ m\geq 1, $  $ k\geq 0,$ and $(m,k)\notin \{(1,0),(1,1),(1,2),(1,3),(2,0),(2,1) \}.$ We consider the following cases.\\
Case $(i).$ If $k=0,$ then $m\geq 3$ and so    $\dim L=2m+2\geq 8.$ Thus $L/N \cong H(m)=\langle  \overline{x_{i}}, \overline{y_{i}}, \overline{z} \mid [\overline{x_{i}}, \overline{y_{i}}]=\overline{z},1\leq i\leq m\rangle, $ in where $ \overline{x_{i}}=x_{i}+N, $ $ \overline{y_{i}}=y_{i}+N $ and $ \overline{z}=z+N $ for all $ i $ such that $ 1\leq i\leq m. $ Put $N=\langle z_{1} \rangle.$ Then the set
$\{x_i, y_i, z, z_{1}|1\leq i \leq m\}$ is a basis of $L$ and
\begin{align*}
&[x_{i}, y_{i}]=z+\beta^{\prime}_{ii} z_{1} ~~\text{for all}~i  ~ \text{such that} ~1\leq i\leq m ~\text{and} ~\beta^{\prime}_{ii}\in \mathbb{F},\\
&[x_{r}, y_{s}]=\alpha_{rs}z_{1} ~~ \text{for all}~r, s~ \text{such that}~r\neq s,~1\leq r,s\leq m  ~\text{and} ~\alpha_{rs}\in \mathbb{F},\\
& [x_{t}, x_{t_{1}}]= \alpha^{\prime}_{tt_{1}}z_{1} ~~ \text{for all}~t,t_{1}~ \text{such that}~1\leq t,t_{1} \leq m ~~\text{and} ~  \alpha^{\prime}_{t,t_{1}} \in \mathbb{F}, \\
& [y_{i_{1}}, y_{j_{1}}]=\beta_{i_{1}j_{1}}z_{1} ~~ \text{for all}~i_1,j_1~ \text{such that}~1\leq i_{1},j_{1} \leq m ~\text{and} ~ \beta_{i_{1}j_{1}}\in \mathbb{F}.
\end{align*}
Put $ z'=z+\beta^{\prime}_{ii} z_{1} .$ A change of a variable allows that   $ \beta^{\prime}_{ii}=0 $ for all $i$ such that $ 1\leq i\leq m $. Since $ \dim L^{2}=2, $ without loss of generality, we may assume that   $ \alpha'_{12}\neq 0. $ Then $  [x_{1}, x_{2}]=\alpha'_{12} z_{1}\neq 0.$ Now a change of variables $ x^{\prime}_{i}= \alpha'^{-1}_{12}x_{i}$, $ y^{\prime}_{i} =\alpha'^{-1}_{12}y_{i}$, $ z'_{1}=\alpha'^{-1}_{12}z_{1} $ and $ z'=\alpha^{-2}_{12}z$ and relabeling, we have
\begin{align*}
&[x_{1}, x_{2}]=z_{1},\cr
&[x_{i}, y_{i}]=z ~~ \text{for all}~i ~\text{such that} ~1\leq i\leq m,\cr
& [x_{t}, x_{t_{1}}]=\alpha^{\prime} _{tt_{1}}z_{1}~~ \text{for all}~t,t_{1} ~\text{such that} ~1\leq t,t_{1} \leq m, ~ (t,t_{1})\notin \lbrace (1,2), (2, 1)\rbrace  \cr
&\text{and} ~\alpha^{\prime} _{tt_{1}}\in \mathbb{F},\cr
&[x_{r}, y_{s}]=\alpha_{rs}z_{1} ~~ \text{for all}~r, s~ \text{such that}~r\neq s,~1\leq r,s\leq m,~~~~  ~\text{and} ~\alpha_{rs}\in \mathbb{F},\cr
& [y_{i_{1}}, y_{j_{1}}]=\beta _{i_{1}j_{1}}z_{1} ~~\text{for all}~i_{1},j_{1} ~\text{such that}~1\leq i_{1},j_{1}\leq m,~ \text{and} ~ \beta _{i_{1}j_{1}}\in \mathbb{F}.
\end{align*}
By  \cite[Lemma 3.3]{N1}, since   $ \dim (L/\langle z\rangle)^{2}=1, $ we get $L/\langle z \rangle \cong H(m^{\prime})\oplus A(2m-2m^{\prime}). $ Let $\overline{x}$ be used to denote  $x+\langle z  \rangle $ for   $x\in L.$ Since  $ \overline{[x_{1}, x_{2}]}=\overline{z_{1} } $, we obtain 
$\overline{[x_{1}, x_{t_{1}}]}=\overline{[x_{1}, y_{s}]} =\overline{[x_{2}, x_{t_{1}}]}= \overline{[x_{2}, y_{s}]} =0$ for all $t_{1},s$ such that $1\leq t_{1},s\leq m $ and  $ t_{1}, s \notin \lbrace 1, 2 \rbrace.$  Hence  $\overline{\alpha^{\prime}_{1t_{1}}z_{1}}=\overline{\alpha_{1s}z_{1}}=\overline{\alpha^{\prime}_{2t_{1}}z_{1}}=\overline{\alpha_{2s}z_{1}}=0$ 
and so $\alpha^{\prime}_{1t_{1}}z_{1}, \alpha_{1s}z_{1}, \alpha^{\prime}_{2t_{1}}z_{1}, \alpha_{2s}z_{1}\in \langle z \rangle.$ Since $ z_{1} $ and $ z $ are elements of the basis, $ \langle z_{1}\rangle \cap \langle z \rangle =0.$ Thus $\alpha^{\prime}_{1t_{1}}=\alpha_{1s}=\alpha^{\prime}_{2t_{1}}=\alpha_{2s}=0$ for all  $s, t_{1}$ such that $1\leq t_{1}, s\leq m $ and  $ t_{1}, s \notin \lbrace 1, 2 \rbrace$. Therefore
$
[x_{1}, x_{t_{1}}]=[x_{1}, y_{s}] =[x_{2}, x_{t_{1}}]= [x_{2}, y_{s}]=0
$
for all $ t_{1},s $ such that $ 1\leq t_{1}, s \leq m $ and  $ t_{1}, s \notin \lbrace 1, 2 \rbrace$.
We claim that $ [x_{r}, y_{s}]=0 $ for all $r,s$ such that $ 3\leq r \leq m $, $ 1\leq s\leq m $ and $ r\neq s $. By contrary, assume that  there is $ \alpha_{r_{1}s_{1}} \neq 0$ for some $ r_{1},s_{1} $ such that $ 3 \leq  r_{1}\leq m $ and $ 2\leq s_{1}\leq m. $ Putting $ x^{\prime}_{r_{1}}=x_{r_{1}}+\alpha_{r_{1}s_{1}} x_{1} $ and $ y^{\prime}_{s_{1}}=y_{s_{1}}-x_{2} $, we obtain that 
\begin{align*}
&[x_{1}, x_{2}]=z_{1},\cr
&[x_{s_{1}},y'_{s_{1}}]=[x'_{r_{1}},y_{r_{1}}]=z, [x_{i}, y_{i}]=z ~~ \text{for all}~i ~\text{such that} ~1\leq i \leq m~ \text{and}~  i\notin \{ r_{1}, s_{1}\},\cr
& [x_{t}, x_{t_{1}}]=\alpha^{\prime} _{tt_{1}}z_{1}~~ \text{for all}~t,t_{1}\neq r_{1} ~\text{such that}~1\leq t,t_{1} \leq m, ~ (t,t_{1})\notin \lbrace (1,2), (2, 1)\rbrace \\ &\text{and} ~\alpha^{\prime} _{tt_{1}}\in \mathbb{F},\cr
& [x_{s}, x'_{r_{1}}]=\alpha' _{sr_{1}}z_{1}, \cr
&[x_{1}, x_{r_{1}}]=\alpha' _{1r_{1}}z_{1},\cr
&[x_{2}, x_{r_{1}}]=\alpha' _{2r_{1}}z_{1},\cr
  & [x_{r},y'_{s_{1}}]=-z_{1},\cr &
   [y'_{i_{1}},y_{j_{1}}]=\beta _{i_{1}j_{1}}z_{1} ~\text{for all}~ \beta _{i_{1}j_{1}}\in \mathbb{F} ~\text{and} ~ i_{1},j_{1} \notin \lbrace s_{1}, 2 \rbrace, \cr
  & [y'_{s_{1}},y_{2}]=\beta_{s_{1}2}z_{1}-z,\cr & [x'_{r_{1}},y_{1}]=\beta z_{1}+\alpha_{r_{1}s_{1}}z.
\end{align*}
Since $ L/\langle z_{1}\rangle \cong H(m)$ and $ \overline{[x_{1},y_{1}]}=\overline{z}, $ we have $\overline{[x^{\prime}_{r_{1}},y_{1}]}=\overline{\alpha_{r_{1}s_{1}} z }=0 $ and so $ \alpha_{r_{1}s_{1}} z \in \langle z_{1}\rangle. $ Since $ z_{1} $ and $ z $ are elements of the basis, $ \langle z_{1}\rangle \cap  \langle z \rangle=\langle 0 \rangle. $
Therefore $ \alpha_{r_{1}s_{1}}=0. $ It is a contradiction. Hence $ [x_{r}, y_{s}]=0 $ for all $r,s$ such that $ 3\leq r \leq m $, $ 2\leq s\leq m $ and $ r\neq s $.\\
If $ s_{1}=1, $  then a change of a variable $ y^{\prime}_{1}=y_{1}-\alpha_{r1}y_{r} $ shows that 
\begin{align*}
  &[x_{r}, y^{\prime}_{1}]=\alpha_{r1}z_{1}-\alpha_{r1}z,\cr
 & [x_{1},y_{1}]=[x_r, y'_{1}]=z,\cr
 &[y'_{1},y_2]=-z.
\end{align*}
Since $ L/ \langle z_{1}\rangle  \cong H(m)$ and $ \overline{[x_{1},y_{1}]}=\overline{z}, $ we have $\overline{[x_{r},y'_{1}]}=\overline{\alpha_{r1} z}=0 $ and  $ \alpha_{r1}=0. $ It is a contradiction, hence $ [x_{r}, y_{1}]=0. $ Therefore $ [x_{i}, y_{j}]=0 $ for all $i,j$ such that $1\leq i,j\leq m$ and  $i\neq j. $
\\ By a similar way, we can see that
$ [y_{i_{1}}, y_{j_{1}}]=[x_{r}, x_{s}]=0 $ for all $i_{1},j_{1},r,s $ such that $ 1\leq i_{1},j_{1},r,s \leq m,$  and $(r,s)\notin \lbrace (1,2),(2, 1)\rbrace.$ Thus
\[[x_{i}, y_{i}]=z ~~\text{for all}~i ~ \text{such that}~ 1\leq i\leq m,~ and ~ [x_{1}, x_{2}]=z_{1}.\]
 Therefore $ L $ is isomorphic to
\[H_{1}=\langle x_{i}, y_{i}, z, z_{1}\mid [x_{i}, y_{i}]=z, [x_{1}, x_{2}]=z_{1}, ~~1\leq i\leq m\rangle. \]
Case $ (ii). $ If $k=1,$ then  $m\geq 3$  and so $\dim L=2m+3> 8.$ Thus $ L/\langle z_{1} \rangle \cong H(m)\oplus A(1)=\langle  \overline{x_{i}}, \overline{y_{i}}, \overline{z} \mid [\overline{x_{i}}, \overline{y_{i}}]=\overline{z},1\leq i\leq m\rangle\oplus A(1).$ Put  $ A(1)\cong \langle q \rangle.$  Hence the set
$\{x_i, y_i, z, z_{1},q|1\leq i \leq m\}$ is a basis of $L$ and
\begin{align*}
&[x_{i}, y_{j}]=\alpha_{ij}z_{1}  ~~\text{for all}~i,j~ \text{such that}~i\neq j, ~1\leq i,j\leq m,  ~~ \text{and}~~ \alpha_{ij}\in \mathbb{F},\cr
& [x_{t}, x_{r}]=\alpha ^{\prime}_{rt}z_{1} ~~\text{for all}~r,t~ \text{such that}~1\leq r,t \leq m ~\text{and} ~ \alpha ^{\prime}_{rt}\in \mathbb{F},\cr
& [y_{i_{1}}, y_{j_{1}}]=\beta_{i_{1}j_{1}}z_{1} ~~\text{for all}~i_{1},j_{1}~ \text{such that}~1\leq i_{1},j_{1} \leq m~ \text{and} ~ \beta_{i_{1}j_{1}}\in \mathbb{F},\cr
&[x_{i}, y_{i}]=z+\beta^{\prime}_{ii} z_{1} ~~\text{for all}~i~  \text{such that}~1\leq i\leq m~\text{and}~ \beta^{\prime}_{ii}\in \mathbb{F}, \cr
& [q, x_{i}]=\alpha_{i} z_{1}~~ \text{for all}~ i ~\text{such that}~ 1\leq i\leq m ~ \text{and} ~ \alpha_{i}\in \mathbb{F},\cr
& [q, y_{i}]=\alpha^{\prime}_{i} z_{1}~~ \text{for all}~i~ \text{such that} ~1\leq i\leq m ~ and ~ \alpha^{\prime}_{i}\in \mathbb{F}.
\end{align*}
 In a similar way as the case $ (i) $, we can see that
 \begin{align*}
 &[y_{i_{1}}, y_{j_{1}}]=0 ~~\text{for all}~i_1,j_1~ \text{such that}~1\leq i_1,j_1\leq m,\cr
 & [x_{t}, x_{r}]=0 ~~\text{for all}~r,t ~\text{such that} ~1\leq r,t \leq m ~\text{and} ~ (r,t)\notin \lbrace (2,1) , (1,2)\rbrace ,\cr
 &[x_{1}, x_{2}]=\alpha_{12}z_{1},\cr
 &[x_{i}, y_{j}]=0~~\text{for all}~ i,j~ \text{such that} ~i\neq j,~\text{and}~1\leq i,j\leq m,  \cr
 & [x_{i}, y_{i}]=z~~\text{for all}~ i~ \text{such that}~1\leq i\leq m, \cr
 & [q, x_{i}]=\alpha_{i} z_{1}~~ \text{for all}~ i ~\text{such that}~ 1\leq i\leq m ~ \text{and} ~ \alpha_{i}\in \mathbb{F},\cr
& [q, y_{i}]=\alpha^{\prime}_{i} z_{1}~~ \text{for all}~i~ \text{such that} ~1\leq i\leq m ~\text{ and} ~ \alpha^{\prime}_{i}\in \mathbb{F}.
 \end{align*}
 Since $q\notin Z(L)=L^2,$  there exists an element $ x_{i} $ or $ y_{i} $  for some $i,~1\leq i\leq m$ such that $ [q,  x_{i}]\neq 0 $ or  $[q, y_{i}]\neq 0.$  Without loss of generality, let $[q,  x_{1}]=\alpha z_{1}\neq 0.$ Put $ x^{\prime}_{i}=\alpha^{-1}x_{i}, $ $ y^{\prime}_{i}=\alpha^{-1}y_{i}, $ $ z'=\alpha^{-2}z $ and $ z'_{1}=\alpha^{-1}z_{1} $ for all $i$ such that $~1\leq i\leq m.$ A change  of variables allows that
 \begin{align*}
 &[x'_{1}, x'_{2}]=\alpha_{12}\alpha^{-1}z_{1},\cr
 & [x^{\prime}_{i}, y^{\prime}_{i}]=z' ~~\text{for all}~ i~~\text{such that}~ 1\leq i\leq m,  \cr
 &[q,x^{\prime}_{1}]=z_{1},\cr
& [q, x'_{i}]=\alpha^{-1} \alpha_{i}z_{1},\cr & [q, y'_{i}]=\alpha^{-1} \alpha'_{i}z_{1} ~~\text{for all}~ i~\text{such that}~1\leq i\leq m.  
 \end{align*}
 By invoking  \cite[Lemma 3.3]{N1} and $\dim L=2m+3,$ we get $ L/\langle z' \rangle \cong H(m^{\prime})\oplus A(2m-2m^{\prime}+1). $  Since  $ \overline{[q, x_{1}]}=\overline{z_{1} } $, we obtain
$\overline{[x_{1}, x_{2}]}=\overline{[q, x_{r}]}=\overline{[q, y_{i}]}=0$ for all  $r,i$ such that $1\leq r,i\leq m $ and so $\alpha_{12}\alpha^{-1}=\alpha_{i}\alpha^{-1}=\alpha'_{i}\alpha^{-1}=0$. Since $\alpha \neq 0,$ we get $\alpha_{12}=\alpha_{i}=\alpha'_{i}=0$ for all $ i $ such that $ 1\leq i\leq m. $ Hence, $[x_{1}, x_{2}]=[q, x_{r}]=[q, y_{i}]=0$ for all  $r,i$ such that $1\leq r,i\leq m $.
 Now by labelling $ z' $ with $ z, $ $ L $ is isomorphic to
 \[H_{2}=\langle x_{i}, y_{i}, q, z, z_{1}\mid [x_{i}, y_{i}]=z, [q, x_{1}]=z_{1}, ~~1\leq i\leq m\rangle. \]
 Case $(iii).$ If $ k\geq 2,$  then $ m\geq 1.$    
 Now let $ m\geq 2, $ we have $\dim L=2m+k+2\geq 8.$ Since $L /\langle z \rangle \cong H(m)\oplus A(k),$  there exist two ideals $ L_{1}/ \langle z \rangle $ and $ L_{2}/ \langle z \rangle $ of $ L/ \langle z \rangle$ where
 \begin{align*}
&  L_{1}/ \langle z \rangle =\langle  \overline{x_{i}}, \overline{y_{i}}, \overline{z_{1}}  \mid [\overline{x_{i}}, \overline{y_{i}}]=\overline{z_{1}},1\leq i\leq m\rangle \cong H(m),\cr
& L_{2}/ \langle z \rangle \cong \bigoplus_{j=1}^k \langle (q_{j}+ \langle z \rangle) \rangle \cong A(k).
 \end{align*}
 Hence the set
$\{x_i, y_i, q_j,z,z_{1}|1\leq j\leq k,1\leq i\leq m\}$ is a basis of $L$ and
  $ L^{2}_{2}\subseteq  \langle z \rangle.$ Thus, we have two subcases.\\
\textbf{Subcase $(a).$ }Since $ L_{2}^{2}=0, $ $ L_{2}=\bigoplus_{j=1}^k \langle q_{j}\rangle\oplus  \langle z \rangle. $ By a similar way used in   cases $ (i) $ and $ (ii), $ we have
\begin{align*}
 &[x_{i}, y_{i}]=z ~~  \text{for all }~i~ \text{such that}~1\leq i\leq m,  \cr
 &[x_{1}, x_{2}]=\alpha z_{1} ~ \text{and} ~ \alpha \in \mathbb{F},\cr
&[x_{i'}, y_{s}]=0~~\text{for all}~i',s~\text{such that}~1\leq i',s\leq m~~ \text{and}~~ i'\neq s, \cr
& [x_{t}, x_{r}]=0~~\text{for all}~r,t~\text{such that}~1\leq r,t \leq m~~\text{and}~~ (t,r)\notin \lbrace (1,2),(2,1)\rbrace ,\cr
& [y_{i_{1}}, y_{s_{1}}]=0 ~~\text{for all}~i_{1},s_{1}~\text{such that}~1\leq i_{1},s_{1} \leq m,\cr
& [q_{j}, x_{s'}]=\alpha_{js'}z_{1} ~~\text{for all}~j,s'~ \text{such that}~ 1\leq s'\leq m, 1\leq j\leq k,~ \text{and} ~ \alpha_{js'}\in \mathbb{F},\cr
& [q_{j}, y_{i}]=\alpha^{\prime}_{ji}z_{1}~~ ~~\text{for all}~i,j~ \text{such that}~ 1\leq i\leq m, 1\leq j\leq k, ~  ~~\text{and}~~\alpha^{\prime}_{ji}\in \mathbb{F}.
\end{align*}
Now let $ [x_{1}, x_{2}]=\alpha z_{1}\neq 0.$ Put $ y^{\prime}_{1}=y_{1}-\alpha'_{s1}x_{1}$ and $ y^{\prime}_{2}=y_{2}-\alpha'_{s2}x_{1}$ for all $s$ such that $ 1\leq s\leq k, $ we have 
\begin{align*}
&[y'_{1}, y'_{2}]=\alpha'_{s2}z, [x_1, y^{\prime}_{1}]=[x_2, y^{\prime}_{2}]=z,
[x_2,y^{\prime}_{1}]=\alpha '_{s1}z,\cr
& [q_{i}, y'_{1}]=\alpha _{i1}z_{1}-\alpha '_{s1}\alpha _{i1}z_{1},    [q_{i}, y'_{2}]=\alpha _{i2}z_{1}-\alpha '_{s2}\alpha _{i2}z_{1}.
\end{align*}
Since $ L/\langle z \rangle \cong H(m')\oplus A(2m+k-2m'+1) $ and $\overline{ [x_{1}, x_{2}]}=\alpha \overline{z_{1}}\neq 0,$  $\alpha_{s1}=\alpha_{s2}=\alpha'_{s1}=\alpha'_{s2}=0 $ for all $s$ such that $ 1\leq s\leq  k.$ Thus $
 [q_{s}, x_1]=[q_{s}, y_{1}]=[q_{s}, x_2]=[q_{s}, y_{2}]=0$ for all $s$ such that $ 1\leq s\leq  k.$
 Since $q_s\notin Z(L)=L^2$ for all $s$ such that  $ 1\leq s\leq  k,$ there exists an element $ x_{i} $ or $ y_{i} $  for some $i$ such that $3\leq i\leq m$ and $ [q_s,  x_{i}]\neq 0 $ or  $[q_s, y_{i}]\neq 0$  for all $s$ such that $ 1\leq s\leq  k.$\\
   We conclude that $2\leq k\leq 2m-4.$ Without loss of generality, let $ [q_{s},x_{s+2}]=\alpha_{s}z_{1}\neq 0$ for all $s$ such that $ 1\leq s\leq  m-2$ and $ [q_{s},y_{s-(m-4)}]=\alpha_{s}z_{1}= 0$ for all  $ m-1\leq s\leq  2m-4.$  A change of a variable allows that $ \alpha_{s}=1$ for all $s$ such that $ 1\leq s\leq  k\leq 2m-4.$ \\
 If $ m-1\leq k\leq 2m-4, $ then $ L $ is isomorphic to 
 \[\langle x_{i}, y_{i}, q_{j}, z, z_{1}\mid [x_{i}, y_{i}]=z, [x_{1}, x_{2}]=[q_{s},x_{s+2}]= [q_{t},y_{t-(m-4)}]=z_{1}, \]\[1\leq i\leq m, 1\leq s\leq m-2, m-1\leq t\leq 2m-4,1\leq j\leq k\rangle.\]
 Put $ q'_{1}=q_{1}-q_{m-1} $ and  $ x'_{3}=x_{3}-y_{3}. $  We have $ [q'_{1}, x'_{3}]=0 ,$ $  [q'_{1}, y_{3}]=-z_{1} $ and  $[q_{m-1}, y_{3}]=z_{1}.$
Since  $ \overline{ [q'_{1},y_{1}]}=\overline{[q_{1}, y_{1}]}=\overline{z_{1}} $ are relations   of  $ L/\langle z \rangle \cong H(m'')\oplus A , $ we have a contradiction  by looking relations of $ L/\langle z \rangle \cong H(m'')\oplus A.  $  
  Therefore $  2\leq k\leq m-2, $ then $ L $ is isomorphic to\\
 \[ H_{3}=\langle x_{i}, y_{i}, q_{j}, z, z_{1}\mid [x_{i}, y_{i}]=z, [x_{1}, x_{2}]=[q_{s},x_{s+2}]=z_{1}, 1\leq i\leq m,\]\[ 1\leq s\leq k, 1\leq j\leq k \rangle ~~\text{for all} ~k,m~~\text{such that}~~2\leq k\leq m-2~~ \text{and}~~ m\geq 4. \]
Let $ [x_{1}, x_{2}]= 0. $ By a similar technique used in cases $ (i) $ and $ (ii) $,
if $ m+1 \leq k\leq 2m, $ then we have a contradiction.
 One can check that for all $ k $ such that $ 2\leq k\leq m, $ $ L $ is isomorphic to 
  \[ H_{4}=\langle x_{i}, y_{i}, q_{j}, z, z_{1}\mid [x_{i}, y_{i}]=z, [q_{s},x_{s}]=z_{1}, 1\leq i\leq m, 1\leq s\leq k,  1\leq j\leq k \rangle\]\[\text{for all}~~k,m~~\text{such that}~~2\leq k \leq m~~\text{and}~~ m\geq 2.\]
\textbf{Subcase $(b).$ }Since $ L_{2}^{2}=\langle z \rangle, $ and $ L_{2}\cong H(k_{1})\oplus A(r), $  we have two subalgebras of $L$ such that
$ L_3=\langle q_{i}, q^{\prime}_{i}, z_{1} \mid  [q_{i},q^{\prime}_{i}]=z_{1}, 1\leq i\leq k_{1} \rangle\cong H(k_{1})$ and $ L_4=\bigoplus \limits_{i=0}^{r} \langle p_{i}\rangle \cong A(r).$ By a similar way used in the subcase $(a),$  one of the following  cases should be occured.\\
\textbf{(b-1).} If $[x_1,x_2]\neq 0$ and $ r=0, $ then  $L$ is isomorphic to
\[  H_{5}=\langle x_{i}, y_{i}, q_{j}, q_{j}^{\prime}, z, z_{1}\mid [x_{i}, y_{i}]=z, [x_{1}, x_{2}]=[q_{j},q^{\prime}_{j}]=z_{1}, 1\leq i\leq m, 1\leq j\leq k_{1} \rangle \]\[ \text{for all}~k_{1}~~\text{such that} ~~k_1\geq 2.\]
\textbf{(b-2).} If $[x_1,x_2]\neq 0$ and $1\leq  r \leq m-2, $ then $L$ is isomorphic to
\[  H_{6}=\langle x_{i}, y_{i}, q_{j}, q_{j}^{\prime}, p_{t},z, z_{1} \mid [x_{i}, y_{i}]=z, [x_{1}, x_{2}]=[q_{j},q^{\prime}_{j}]=z_{1}, \]\[[p_{s},x_{s+2}]=z, 1\leq i\leq m, 1\leq j\leq k_{1}, 1\leq t\leq r \rangle. \]\\
\textbf{(b-3).} If $[x_1,x_2]= 0$ and $ r=0, $ then  $L$ is isomorphic to
\begin{align}\label{e7}
&\langle x_{i}, y_{i}, q_{j}, q_{j}^{\prime}, z, z_{1}\mid [x_{i}, y_{i}]=z, [q_{j},q^{\prime}_{j}]=z_{1}, 1\leq i\leq m, 1\leq j\leq k_{1}\rangle \cr
&\cong H(m)\oplus H(k_{1}). 
\end{align}
 \textbf{(b-4).} If $[x_1,x_2]= 0$ and $1\leq  r\leq m,$  then $L$ is isomorphic to
\[ H_{7}=\langle x_{i}, y_{i}, q_{j}, q_{j}^{\prime}, p_{t},z, z_{1}\mid [x_{i}, y_{i}]=z,[q_{j},q^{\prime}_{j}]=[p_{t},x_{t}]=z_{1}, \]\[1\leq i\leq m, 1\leq j\leq  k_{1}, 1\leq t\leq r\rangle. \]
Now let  $ m=1.$  Then $ i=1 $ and $ k\geq 4. $ If $ L_{2} $ is abelian, then $ [q_{j}, q'_{j}]=0 $ for all $ j $ such that  $ 1\leq  j \leq k. $  Since $ \dim L\geq 8, $ we have  $ j\geq 4. $ On the other hand  since $q_j\notin Z(L)=L^2$ for all $j$ such that  $ 1\leq j \leq k,$ there exists an element $ x_{1} $ or $ y_{1} $  such that  $ [q_j,  x_{1}]\neq 0 $ or  $[q_j, y_{1}]\neq 0$  for all $j$ such that $ 1\leq j \leq k.$ 
 Hence there is $ j $ for some $ 1 \leq j\leq k $  such that  $ q_{j} \in  Z(L).$ It is a contradiction. Thus this case does not happen.\\
Let $ \dim L_{2}^{2}=1, $ then $ L_{2}\cong H(k_1)\oplus A(r). $ \\
Consider $ L_3=\langle q_{j}, q^{\prime}_{j}, z_{1} \mid  [q_{j},q^{\prime}_{j}]=z_{1}, 1\leq j\leq k_{1} \rangle\cong H(k_{1}),$ $ L_4=\bigoplus \limits_{s=0}^{r} \langle p_{s}\rangle \cong A(r) $ and  $ 0\leq r \leq m.$ Since  $m=1, $ we have
$ 0\leq r \leq 1. $\\
If $ r=0, $  then $ L $ is isomorphic to 
\begin{align}\label{e9}
\langle x_{1}, y_{1}, q_{i}, q'_{i}, z, z_{1}\mid [x_{1}, y_{1}]=z, [q_{i},q^{\prime}_{i}]=z_{1}, 1\leq i\leq k_{1}\rangle\cong H(1)\oplus H(k_{1}) 
\end{align}
for all $ k_{1} $ such that $ k_{1}\geq 2$.\\
The structures of (\ref{e7}) and (\ref{e9}) imply
\[  H_{8}=\langle x_{i}, y_{i}, q_{j}, q_{j}^{\prime}, z, z_{1}\mid [x_{i}, y_{i}]=z, [q_{j},q^{\prime}_{j}]=z_{1}, 1\leq i\leq m, 1\leq j\leq k_{1}\rangle\]\[\cong H(m)\oplus H(k_{1})~~\text{for all}~~ k_{1}, m~~\text{such that} ~~ k_{1}\geq 2~~\text{and} ~~m\geq 1.\]
If $ r=1, $ then $ L $ is isomorphic to 
\begin{align}\label{e3}
\langle x_{1}, y_{1}, q_{i}, q_{i}^{\prime}, p_{1}, z, z_{1}\mid [x_{1}, y_{1}]=z, [q_{i},q^{\prime}_{i}]=[p_{1},x_{1}]=z_{1}, 1\leq i\leq k_{1} \rangle.
\end{align}
By considering a homomorphism  that maps $ x_{1}\mapsto x_{1}, $ $ y_{1}\mapsto y_{1}, $ $ p_{1}\mapsto q, $ 
$ q_{i}\mapsto x_{i}, $ and $ q'_{i}\mapsto y_{i} $ for all $ i $ such that $ i\geq 4, $ the structures of 
(\ref{e3}) and $ H_{2} $ are isomorphic.\\
The result follows.
\end{proof}
\end{thm}
\begin{cor}
Let $H$ be a  generalized Heisenberg Lie algebra  of rank $ 2 $ over an arbitrary field. Then  $H $ is isomorphic to one of the Lie algebras 
 $ L_{5,8}, $ $ L_{6,22}(\varepsilon), $ $ L^{(2)}_{6,7}(\eta), $ $ L_{1}, $ $ L_{2}, $ or $ H_{i} $ for all $ i $ such that $ 1\leq i\leq 8. $
 \begin{proof}
Using Theorem \ref{T1} and \cite[Theorems 2.7 and 3.1]{N3}, the result  follows.
 \end{proof}
\end{cor}
For a  Lie algebra $ L $, the concept epicenter $  Z^{*}(L) $ was introduced in \cite{sal}. The importance of  $  Z^{*}(L) $ is due to the fact that $ L $ is capable if and only if $  Z^{*}(L)=0 $. Another notion having relation  to capability is the exterior center $  Z^{\wedge}(L) $ of a Lie algebra $ L. $ The exterior center $ Z^{\wedge}(L) $ is the set of all elements $ l $ of $ L $ for which $ l\wedge l'=0 $ for all $ l'\in L .$  It is known that from   \cite[Lemma 3.1]{N5}, $ Z^{*}(L)=Z^{\wedge}(L) $. A Lie algebra $ L $ is called unicenrtal if $ Z^{\wedge}(L)=Z(L). $ \\
In the following theorem, we determine the exterior center of $ H_{i} $ for all $ i $ such that  $ 1\leq i\leq 8. $
\begin{thm}
 Let $H$ be an $ n $-dimensional  generalized Heisenberg Lie algebra  of rank $ 2 $ over an arbitrary field $ \mathbb{F} $ for all $ n $ such that  $ n\geq 8. $ Then 
\begin{itemize}
\item[$(i)$]$H$ is unicentral if and only if $H$ is isomorphic to one of the  Lie algebras  $H_3,H_4,H_5,H_6,$  $ H_{7}, $ or $H_8=H(m)\oplus H(k_{1}) $ for all $ k_{1}, m $ such that  $   k_{1}\geq 2  $  and $ m\geq 2, $
\item[$(ii)$]$Z^{\wedge}(H)\cong A(1) $ if and only if $H$ is isomorphic to one of the  Lie algebras   $H_1,$  $H_2,$ or $ H_{8}=H(1)\oplus H(k_{1}) $ for all $ k_{1} $ such that $ k_{1}\geq 2. $
\end{itemize}
\begin{proof}
  Assume that $ H$ is isomorphic to 
  \[ H_{3}=\langle x_{i}, y_{i}, z, z_{1}, q_{j}\mid [x_{i}, y_{i}]=z, [x_{1}, x_{2}]=[q_{j}, x_{j+2}]=z_{1}, ~~1\leq i\leq m, 1\leq j \leq k\rangle \]
  for all $ k,m $ such that  $2\leq k\leq m-2$ and $m\geq 3.$  We claim that $Z^{\wedge}(H_3)=H_3^2.$ For all $ i,j, $ such that $i\neq j, 1\leq i,j\leq m,$ and $(i,j)\neq (1,2), $ we have
  \begin{align*}
 & x_i\wedge z =x_i\wedge [x_j,y_j]=[x_i,y_j]\wedge x_j- [x_i,x_j]\wedge y_j=0,\cr
  & y_i\wedge z=y_i\wedge [x_j,y_j]=[y_i,y_j]\wedge x_j- [y_i,x_j]\wedge y_j=0,\cr
    &q_j\wedge z =q_j\wedge [x_1,y_1]=[q_j,y_1]\wedge x_1- [q_j,x_1]\wedge y_1=0~~ \text{for all }~~j~~\text{such that}~~1\leq j\leq k,
  \end{align*}
  and 
  \begin{align*}
 & x_1\wedge z_{1}=x_1\wedge [q_1,x_3]=[x_1,q_1]\wedge x_3- [x_1,x_3]\wedge q_1=0,\cr
  & x_2\wedge  z_{1}=x_2\wedge [q_1,x_3]=[x_2,q_1]\wedge x_3- [x_2,x_3]\wedge q_1=0,\cr
   &y_1\wedge  z_{1}=y_1\wedge [q_1,x_3]=[y_1,q_1]\wedge x_3- [y_1,x_3]\wedge q_1=0,\cr
  &y_2\wedge  z_{1}=y_2\wedge [q_1,x_3]=[y_2,q_1]\wedge x_3- [y_2,x_3]\wedge q_1=0,\cr
   &x_t\wedge  z_{1}=[x_t,x_1]\wedge x_2- [x_t,x_2]\wedge x_1=0~~\text{for all }~~t~~\text{such that}~~3\leq t\leq m,\cr
    &y_t\wedge  z_{1}=[y_t,x_1]\wedge x_2- [y_t,x_2]\wedge x_1=0~~\text{for all }~~t~~\text{such that}~~3\leq t\leq m,\cr
  & q_j\wedge  z_{1}=[q_j,x_1]\wedge x_2- [q_j,x_2]\wedge x_1=0~~ \text{for all }~~r~~\text{such that}~~1\leq j\leq k,\cr
    &z\wedge  z_{1}=0.
  \end{align*}
  Hence $ z $ and $  z_{1} $ are in $ Z^{\wedge}(H_3) $ and so $Z^{\wedge}(H_3)=Z(H_3)=H_3^2.$ Therefore    $H_3$ is unicentral. By a similar way, $H_3,H_4,H_5,H_6,$ $ H_{7}, $ and $H_8=H(m)\oplus H(k_{1}) $ for all $ k_{1},m $ such that $  k_{1}\geq 2 $ and $ m\geq 2 $  are unicentral.
  Also, we can obtain that $Z^{\wedge}(H_8)=Z^{\wedge}(H(1)\oplus H(k_{1})) =\langle z \rangle $ for all $ k_{1} $ such that $ k_{1}\geq 2$  and   $Z^{\wedge}(H_1)=Z^{\wedge}(H_2)=\langle z \rangle.$ The proof is completed. 
\end{proof}
\end{thm}
The Schur multiplier of a Lie algebra $ L $ was defined as $ \mathcal{M}(L)\cong ( R\cap F^{2})/ [R,F] $ such that $ L\cong F/R $ for a free Lie algebra $ F. $ The reader can find more information  in \cite{S,H,N1,sal}.
Following corollaries and \cite[Theorem 1]{S} enable us to compute the Schur multiplier of nilpotent Lie algebras of class two with derived subalgebra of dimension two.\\
Let $ H $ be a generalized Heisenberg Lie algebra  of rank $ 2.$ Then Proposition $3.1$ of \cite{N1} shows that $\dim \mathcal{M}(H)=\frac{1}{2}(n-3)(n-2)-2$ when $ Z^{\wedge}(H)=H^{2} $ and $\dim \mathcal{M}(H)=\frac{1}{2}(n-1)(n-4)+1$ when $ Z^{\wedge}(H)\cong A(1). $  In the following corollary, the Schur multiplier of  $H_{i}$ for all $ i $ such that $ 1\leq i\leq 8 $ are computed. 
\begin{cor}
Let $H_{i}$ be an $ n $-dimensional  generalized Heisenberg Lie algebra  of rank $ 2 $ over an arbitrary field  $ \mathbb{F} $ for all $ n,i $ such that $ n\geq 8 $ and $ 1\leq i\leq 8. $ Then
\begin{itemize}
\item[(i).] $  \mathcal{M}(H_{1})\cong \mathcal{M}(H_{2})\cong  A(\frac{1}{2}(n-1)(n-4)+1),$
\item[(ii).]$\mathcal{M}(H_{8})\cong \mathcal{M}(H(1)\oplus H(k))\cong  A(\frac{1}{2}(n-1)(n-4)+1) $ for all $ k $ such that $ k\geq 2,$
\item[(iii).]$ \mathcal{M}(H_{i})\cong A(\frac{1}{2}(n-3)(n-2)-2)$ for all $ i $ such that $3\leq i\leq 7,$
\item[(iv).]$\mathcal{M}(H_{8})\cong \mathcal{M}(H(m)\oplus H(k_1))\cong  A(\frac{1}{2}(n-3)(n-2)-2)$ for all $ k_{1},m $ such that $ k_1\geq 2$  and $ m\geq 2. $ 
\end{itemize}
\begin{proof}
The proof is similar to \cite[Proposition 3.1]{N1}.
\end{proof}
\end{cor}
In the following corollary, we obtain the structure of  all nilpotent  Lie algebras of class two over an arbitrary field with derived subalgebra of dimension two.
 \begin{cor}
Let $L$  be an $ n $-dimensional  nilpotent  Lie algebra of class two over an arbitrary field $ \mathbb{F} $ with derived subalgebra of dimension two. Then $ L $ is isomorphic to one of the following Lie algebras.
\begin{itemize}
\item[(i).] $ L_{5,8}\oplus A(n-5),$
\item[(ii).] $ L_{6,22}(\varepsilon)\oplus A(n-6),$
\item[(iii).]$ L^{(2)}_{6,7}(\eta) \oplus A(n-6),$
\item[(iv).] $ L_{1} \oplus A(n-7),$ 
\item[(v).] $ L_{2} \oplus A(n-7),$
\item[(vi).] $H_{i}\oplus A(n-t)$ for all $ i, t $ such that $ 1\leq i \leq 8, $ and $ t\geq 8. $
\end{itemize}
\begin{proof}
The result follows  from \cite[Theorems 2.7 and 3.1]{N3}, Theorem \ref{T1} and   \cite[Proposition 2.4 ]{N2}.
\end{proof}
\end{cor}

\section*{Acknowledgment}
The  first  author  was  supported  by postdoctoral grant  `` CAPES/PRINT-Edital $n^\circ 41/2017$, Process number:88887.364925/2019-00, ''  at  Federal University of Minas Gerais.

\end{document}